\documentclass{amsart}
\usepackage{amsfonts}
\usepackage{amsmath}
\usepackage{graphicx}
\usepackage[english]{babel}

\newcommand{\bohrequiv}{Bohr-equivalen}
\newcommand{\newequiv}{$^*$-equivalen}

\newtheorem{theorem}{Theorem}

\newtheorem{conjecture}[theorem]{Conjecture}
\newtheorem{corollary}[theorem]{Corollary}

\newtheorem{definition}[theorem]{Definition}

\newtheorem{lemma}[theorem]{Lemma}

\newtheorem*{theorem*}{Theorem}
\newtheorem*{lemma*}{Lemma}

\theoremstyle{definition}
\newtheorem*{remark*}{Remark}

\begin{document}

\title{The converse of Bohr's equivalence theorem with Fourier exponents linearly independent over the rational numbers}

\author{M. Righetti}
\address{Department of Mathematics\\
University of
Genova,
Genoa\\
Italy} \email{righetti@dima.unige.it}

\author{J.M. Sepulcre}
\address{Department of Mathematics\\ University of
Alicante, 03080-Alicante\\
Spain} \email{JM.Sepulcre@ua.es}

\author{T. Vidal}
\address{University of
Alicante, 03080-Alicante\\
Spain} \email{tmvg@alu.ua.es}

\subjclass[2010]{Primary: 42A75, 30D20, 11J72, 11K60}

\keywords{Bohr equivalence theorem; Dirichlet series; Converse theorem; Almost periodic functions}


\begin{abstract}
Given two arbitrary almost periodic functions with associated Fourier exponents which are linearly independent over the rational numbers, we prove that the existence of a common open vertical strip $V$, where both functions assume the same set of values on every open vertical substrip included in $V$, is a necessary and sufficient condition for both functions to have the same region of almost periodicity and to be $^*$-equivalent or Bohr-equivalent. This result represents the converse of Bohr's equivalence theorem for this particular case.\end{abstract}

\maketitle

\section{Introduction}

The theory of almost periodic functions with complex values, created by H. Bohr during the 1920's, opened a way to study a wide class of trigonometric series of the general type and even exponential series. This subject, widely treated in several monographs, has been developed by many authors and has had noteworthy applications \cite{Amerio,Besi,Bochner,Bohr,Corduneanu1,Jessen,Sepulcre-Vidal0}.

The space of almost periodic functions in a vertical strip $U=\{s=\sigma+it:\alpha<\sigma<\beta\}$, $-\infty\leq\alpha<\beta\leq\infty$, which will be denoted in this paper as $AP(U,\mathbb{C})$, is defined as the set of analytic functions $f:U\mapsto \mathbb{C}$ that are equipped with a relatively dense set of almost periods (as Bohr called them) in the following sense: for any
 $\varepsilon>0$ and every reduced strip $U_1=\{s=\sigma+it:\sigma_1\leq\sigma\leq\sigma_2\}$ of $U$ there exists a number $l=l(\varepsilon)>0$ such that every interval of length $l$ contains a number $\tau$ satisfying the inequality
$|f(s+i\tau)-f(s)|\leq \varepsilon$ for all $s$ in $U_1$. In an equivalent way, the space $AP(U,\mathbb{C})$ coincides with the completion of the space of all finite exponential sums of the form $$a_1e^{\lambda_1s}+a_2e^{\lambda_2s}+\ldots+a_ne^{\lambda_ns},$$ with complex coefficients $a_j$ and real exponents $\lambda_j$, equipped with the norm of uniform convergence on every reduced strip of $U$ \cite[p. 148]{Besi}.

Taking as starting point the mean value theorem, the theory of Fourier expansions of periodic functions can be extended to almost periodic functions. Indeed, every function in $AP(U,\mathbb{C})$ can be associated with a certain exponential series of the form $\sum_{n\geq 1} a_ne^{\lambda_n s}$, with complex coefficients $a_n$ and real exponents $\lambda_n$ (the Fourier exponents), which is called the Dirichlet series of the given almost periodic function (see \cite[p.147]{Besi}, \cite[p.77]{Corduneanu1} or \cite[p.312]{Jessen}), and the restriction of this series to vertical lines provides the Fourier series of this function.

In this context, we recall that the class of general Dirichlet series consists of series that take the form
$\sum_{n\geq 1}a_ne^{-\lambda_n s},\ a_n\in\mathbb{C},$
 where
$\{\lambda_n\}$ is a strictly increasing sequence of positive numbers
tending to infinity. Regarding these series, H. Bohr introduced an equivalence relation (which we will refer to as {\bohrequiv}ce) among them that led to exceptional results such as Bohr's equivalence theorem: {\bohrequiv}t general Dirichlet series take the same values in certain vertical lines or strips in the complex plane (see for example \cite{Apostol}). This equivalence relation was used by Righetti in 2017 to obtain a partial converse theorem for the case of general Dirichlet series in their half-plane of absolute convergence \cite{Righetti}.

Regarding the so-called Dirichlet series associated with an almost periodic function $f(s)$ in $AP(U,\mathbb{C})$, it is worth mentioning that $f(s)$ coincides with its associated Dirichlet series in the case of uniform convergence on its strip of almost periodicity $U$ (hence in particular if the convergence is absolute). However, if this condition is not satisfied, we only can state that $f(s)$ is associated with its Dirichlet series on the region $U$. In fact, these Dirichlet series may not converge in $U$ with the ordinary summation, but there exists another way of summation, called Bochner-Fej\'{e}r procedure, which gives rise to a sequence of finite exponential sums, connected with the Dirichlet series, that converges uniformly to $f$ in every reduced strip in $U$, and converges formally to the Dirichlet series on $U$ \cite[p. 148]{Besi}.

More generally, concerning exponential sums of type
$$a_1e^{\lambda_1s}+a_2e^{\lambda_2s}+\ldots+a_je^{\lambda_js}+\ldots,$$
with $a_j\in\mathbb{C}$ and $\{\lambda_1,\lambda_2,\ldots,\lambda_j,\ldots\}$ an arbitrary countable set of distinct real numbers (not necessarily unbounded),
Sepulcre and Vidal established in 2018 a new equivalence relation on them (that we will call {\newequiv}ce, see definitions \ref{DefEquiv00} and \ref{DefEquiv2}), and they also extended it to the context of the complex functions which can be represented by a Dirichlet-like series (in particular those almost periodic functions in $AP(U,\mathbb{C})$) in order to obtain a refined characterization of almost periodicity (see \cite[Theorem 5]{Sepulcre-Vidal0}). This development also led them to an extension of Bohr's equivalence theorem to the case of functions in $AP(U,\mathbb{C})$, which is valid in every open half-plane or open vertical strip included in their region of almost periodicity (under the assumption of existence of an integral basis \cite[Theorem 1]{Sepulcre-Vidal1} and in the general case \cite[Theorem 1]{Sepulcre-Vidal2}). It is convenient to remark that this new {\newequiv}ce relation, which can be formally applied to every Dirichlet series associated with almost periodic functions, coincides with {\bohrequiv}ce \cite{Apostol} (and hence that used in \cite{Righetti}) for the particular case of general Dirichlet series whose sets of exponents have an integral basis.

Given two arbitrary almost periodic functions with associated Fourier exponents which are linearly independent over the rational numbers, the main result in this paper states that they are {\newequiv}t (or also Bohr-equivalent) if and only if there exists an open vertical strip $V$, included in their common region of almost periodicity, where both functions assume the same set of values on every open vertical substrip included in $V$ (see theorems \ref{mth} and \ref{th12}). Also,  we extend this result to the possibility that one of the Fourier exponents is equal to $0$ (see Theorem \ref{mth2}). In fact, we prove that the existence of such an open vertical strip is a necessary and sufficient condition for both functions to have the same region of almost periodicity and to be {\newequiv}t.

Despite the fact that the converse of Bohr's equivalence theorem is, in general, false (see e.g. \cite{Righetti}),
our main result shows that it is true under these conditions on the Fourier exponents (also for the converse of \cite[Theorem 1]{Sepulcre-Vidal2}). In fact, our main theorem is \textit{stronger} than a converse of Bohr's equivalence theorem for this case because it is not necessary to have the same set of exponents.

\section{Preliminaries}

We first consider the following equivalence relation which constitutes our starting point.
\begin{definition}[{\bohrequiv}ce]\label{DefEquiv0}
Let $\Lambda$ be an arbitrary countable subset of distinct real numbers, $V$ the $\mathbb{Q}$-vector space generated by $\Lambda$ ($V\subset \mathbb{R}$), and $\mathcal{F}$
the $\mathbb{C}$-vector space of arbitrary functions $\Lambda\to\mathbb{C}$. 
We define
a relation $\sim$ on $\mathcal{F}$ by $a\sim b$ if there exists a $\mathbb{Q}$-linear map $\psi:V\to\mathbb{R}$ such that
$$b(\lambda)=a(\lambda)e^{i\psi(\lambda)},\ \mbox{with }\lambda\in\Lambda.$$
\end{definition}

The reader may check that this equivalence relation is based on that of Bohr for general Dirichlet series (see e.g. \cite[p. 173]{Apostol}).

Now, let $\Lambda=\{\lambda_1,\lambda_2,\ldots,\lambda_j,\ldots\}$ be an arbitrary countable set of distinct real numbers. We will handle formal exponential sums of the type
 \begin{equation}\label{eqqnew}
\sum_{j\geq 1}a_je^{\lambda_js},\ a_j\in\mathbb{C},\ \lambda_j\in\Lambda,
\end{equation}
where $s=\sigma+it\in\mathbb{C}$. In this context, we will say that $\Lambda$ is a set of exponents and $a_1,a_2,\ldots,a_j,\ldots$ are the coefficients of this exponential sum.

In this way, based on Definition \ref{DefEquiv0}, we consider the following equivalence relation on the classes of exponential sums of type \eqref{eqqnew}.
We will denote as $\sharp \Lambda$ the cardinal of the numerable set $\Lambda$.

\begin{definition}[{\newequiv}ce for exponential sums]\label{DefEquiv00}
Given an arbitrary countable set of distinct real numbers $\Lambda=\{\lambda_1,\lambda_2,\ldots,\lambda_j,\ldots\}$, consider $A_1(s)$ and $A_2(s)$ two exponential sums of the type
$\sum_{j\geq1}a_je^{i\lambda_js}$ and $\sum_{j\geq1}b_je^{i\lambda_js}$, respectively.
We will say that $A_1$ is {\newequiv}t to $A_2$ (in that case, we will write $A_1\ \shortstack{$_{{\fontsize{6}{7}\selectfont *}}$\\$\sim$}\ A_2$) if for each integer value $n\geq 1$, with $n\leq \sharp\Lambda$,
there exists a $\mathbb{Q}$-linear map $\psi_n:V_n\to\mathbb{R}$, where $V_n$ is the $\mathbb{Q}$-vector space generated by $\{\lambda_1,\lambda_2,\ldots,\lambda_n\}$, such that
$$b_j=a_je^{i\psi_n(\lambda_j)},\ j=1,\ldots,n.$$
\end{definition}

Note that this equivalence relation was already introduced in \cite{Sepulcre-Vidal0}, \cite{Sepulcre-Vidal00} and \cite{Sepulcre-Vidal1}. As it was showed in \cite[Proposition 1]{Sepulcre-Vidal00}, it can be characterized in terms of
a basis of the $\mathbb{Q}$-vector space generated by a set $\Lambda=\{\lambda_1,\lambda_2,\ldots\}$ of exponents.
If $G_{\Lambda}=\{g_1,g_2,\ldots\}$ is such a basis, then each $\lambda_j$ in $\Lambda$ is expressible as a finite linear combination of terms of $G_{\Lambda}$, say
\begin{equation*}\label{errej}
\lambda_j=\sum_{k=1}^{i_j}r_{j,k}g_k,\ \mbox{for some }r_{j,k}\in\mathbb{Q},\ i_j\in\mathbb{N},
\end{equation*}
and it is said that $G_{\Lambda}$ is an integral basis for $\Lambda$ if $r_{j,k}\in\mathbb{Z}$ for each $j,k$.

Although definitions \ref{DefEquiv0} and \ref{DefEquiv00} are not equivalent in the general case, it is worth noting that they
are equivalent when it is feasible to obtain an integral basis for the set of exponents $\Lambda$ (see \cite[Proposition 1]{Sepulcre-Vidal2}). For example, this equivalence happens particularly when all the exponents are linearly independent over the rational numbers.

\smallskip

Now we extend Definition \ref{DefEquiv00} to the case of the almost periodic functions in the classes $AP(U,\mathbb{C})$.

\begin{definition}[{\newequiv}ce for almost periodic functions]\label{DefEquiv2}
Given $\Lambda=\{\lambda_1,\lambda_2,\ldots,\lambda_j,\ldots\}$ a set of exponents, let $f_1$ and $f_2$ denote two  functions in $AP(U,\mathbb{C})$, with $U=\{s=\sigma+it:\alpha<\sigma<\beta\}$, whose Dirichlet series are respectively given by
 \begin{equation*}\label{eqq00}
A_1(s)=\sum_{j\geq 1}a_je^{\lambda_js}\ \mbox{and}\ A_2(s)=\sum_{j\geq 1}b_je^{\lambda_js},\ a_j,b_j\in\mathbb{C},\ \lambda_j\in\Lambda.
\end{equation*}
We will say that $f_1$ is {\newequiv}t to $f_2$ if $A_1\ \shortstack{$_{{\fontsize{6}{7}\selectfont *}}$\\$\sim$}\  A_2$, where  $ \shortstack{$_{{\fontsize{6}{7}\selectfont *}}$\\$\sim$}$ is as in Definition \ref{DefEquiv00}. In this case we also write  $f_1\ \shortstack{$_{{\fontsize{6}{7}\selectfont *}}$\\$\sim$}\  f_2$.
\end{definition}

As one can see, the {\newequiv}ce of formal exponential sums (Definition \ref{DefEquiv00}) is the same as the above one for Dirichlet series of almost periodic functions in $AP(U,\mathbb{C})$; this is why it makes sense to use the same notation. 
More generally, {\newequiv}ce can be adapted to the case of the functions (or classes of functions) which are identifiable by their also called Dirichlet series (see \cite[Definition 5]{Sepulcre-VidalUlt} or \cite[Definition 5]{Sepulcre-Vidal00} referred to Besicovitch spaces).

If $f_1$ and $f_2$ are two {\newequiv}t almost periodic functions in $AP(U,\mathbb{C})$, with $U=\{\sigma+it\in\mathbb{C}:\alpha<\sigma<\beta\}$, and $E$ is an open subset of $(\alpha,\beta)$, we recall that, in the same terms as Bohr's equivalence theorem,  the result \cite[Theorem 1]{Sepulcre-Vidal2} assures that the functions $f_1$ and $f_2$ have the same set of values on the region $\{s=\sigma+it\in\mathbb{C}:\sigma\in E\}$. We will deal with the converse of this result for a particular class of functions in $AP(U,\mathbb{C})$.

\section{The closure of the set of values of almost periodic functions}

Given a complex function $f(s)$ and $\sigma_0\in\mathbb{R}$, take the notation $$\operatorname{Img}\left(f(\sigma_0+it)\right)=\{s\in\mathbb{C}:\exists t\in\mathbb{R}\mbox{ such that }s=f(\sigma_0+it)\}.$$

Let $f_1,f_2\in AP(U,\mathbb{C})$ be two {\newequiv}t almost periodic functions in a common vertical strip $U=\{\sigma+it\in\mathbb{C}:\alpha<\sigma<\beta\}$. If $\alpha<\sigma_0<\beta$, we know by \cite[Proposition 4, i)]{Sepulcre-Vidal2} that
$$\overline{\operatorname{Img}\left(f_1(\sigma_0+it)\right)}= \overline{\operatorname{Img}\left(f_2(\sigma_0+it)\right)}.$$
In this section, we will study the validity of this equality for every $\sigma_0\in(\alpha,\beta)$ in terms of the set of values which take $f_1$ and $f_2$ on every region of the form $\{s=\sigma+it\in\mathbb{C}:\sigma\in E\}$, with $E$ an open set of real numbers included in $(\alpha,\beta)$.

\begin{lemma}\label{DeEgea}
Let $f\in AP(U,\mathbb{C})$ with $U=\{\sigma+it\in\mathbb{C}:\alpha<\sigma<\beta\}$, and take $\sigma_0$ such that $\alpha<\sigma_0<\beta$.
Then a complex number $w$ is in $\overline{\operatorname{Img}\left(f(\sigma_0+it)\right)}$ if and only if there exists $\varepsilon_0>0$ satisfying
$$w\in\bigcup_{\sigma\in E_{\sigma_0,\varepsilon}}\operatorname{Img}\left(f(\sigma+it)\right)\ \mbox{for every }0<\varepsilon<\varepsilon_0,$$
where $E_{\sigma_0,\varepsilon}=(\sigma_0-\varepsilon,\sigma_0+\varepsilon)$.
\end{lemma}
\begin{proof}
Let $w_0\in\overline{\operatorname{Img}\left(f(\sigma_0+it)\right)}$, which yields the existence of a sequence $\{t_n\}$ of real numbers such that
$$w_0=\lim_{n\to\infty}f(\sigma_0+it_n).$$ Given $n\in\mathbb{N}$, take the function $h_n(s):=f(s+it_n)$, $s\in U$.
By \cite[Proposition 4]{Sepulcre-Vidal0}, there exists a subsequence $\{h_{n_k}\}_k\subset \{h_n\}_n$ which converges uniformly on reduced strips of $U$ to a function $h(s)$, with $h\ \shortstack{$_{{\fontsize{6}{7}\selectfont *}}$\\$\sim$}\ f$. Note that $$\lim_{k\to\infty}h_{n_k}(\sigma_0)=h(\sigma_0)=w_0.$$ Therefore, by Hurwitz's theorem, there is a positive integer $k_0$ such that for $k>k_0$ the functions $h^*_{n_k}(s):=h_{n_k}(s)-w_0$ have at least one zero in $D(\sigma_0,\varepsilon)$ for every $\varepsilon>0$ sufficiently small. This means that for $k>k_0$ the functions $h_{n_k}(s)=f(s+it_{n_k})$, and hence the function $f(s)$, take the value $w_0$ on the region $\{s=\sigma+it:\sigma_0-\varepsilon<\sigma<\sigma_0+\varepsilon\}$ for every $\varepsilon>0$ sufficiently small.
Consequently, there exists $\varepsilon_0>0$ such that $w_0\in \bigcup_{\sigma\in E_{\sigma_0,\varepsilon}}\operatorname{Img}\left(f(\sigma+it)\right)$, where $E_{\sigma_0,\varepsilon}=(\sigma_0-\varepsilon,\sigma_0+\varepsilon)$ with $0<\varepsilon<\varepsilon_0$ ($\varepsilon>0$ is chosen so that $E_{\sigma_0,\varepsilon}\subset (\alpha,\beta)$).
\smallskip

\noindent Conversely, suppose that $w_0\in\bigcup_{\sigma\in (\sigma_0-\varepsilon,\sigma_0+\varepsilon)}\operatorname{Img}\left(f(\sigma+it)\right)$ for every $0<\varepsilon<\varepsilon_0$ (with $\alpha<\sigma_0-\varepsilon_0<\sigma_0+\varepsilon_0<\beta$). In this way, for each integer value of $n\geq n_0$ with $n_0$ sufficiently large, we have $w_0=f(s_n)$ for some $s_n=\sigma_n+it_n$, with $\sigma_0-\frac{1}{n}<\sigma_n<\sigma_0+\frac{1}{n}$.
Now, let $M$ be an upper bound for $|f'(s)|$ in the region $\{\sigma+it\in\mathbb{C}:\sigma_0-\varepsilon_0<\sigma<\sigma_0+\varepsilon_0\}$ (note that $f'(s)$ is also almost periodic and hence it is bounded on this region \cite[p. 142-144]{Besi}). Therefore, if $n\geq n_0$, we have that
$$|w_0-f(\sigma_0+it_n)|=|f(\sigma_n+it_n)-f(\sigma_0+it_n)|=\left|\int_{\sigma_0}^{\sigma_n}f'(\sigma+it_n)d\sigma\right|\leq M|\sigma_n-\sigma_0|\leq\frac{ M}{ n}.$$
This means that $\lim_{n\to\infty}f(\sigma_0+it_n)=w_0$ and, consequently, $w_0\in\overline{\operatorname{Img}\left(f(\sigma_0+it)\right)}$.
\end{proof}

\begin{theorem}[Equality of the closures of the set of values of almost periodic functions]\label{messi}\
Let $f_1\in AP(U_1,\mathbb{C})$ and $f_2\in AP(U_2,\mathbb{C})$, with $U_1=\{\sigma+it\in\mathbb{C}:\alpha_1<\sigma<\beta_1\}$ and $U_2=\{\sigma+it\in\mathbb{C}:\alpha_2<\sigma<\beta_2\}$ such that $U_1\cap U_2\neq \emptyset$. Consider an interval $(\alpha,\beta)\subset(\alpha_1,\beta_1)\cap (\alpha_2,\beta_2)$.
Then the functions $f_1$ and $f_2$ take the same set of values on every region $\{s=\sigma+it\in\mathbb{C}:\sigma\in E\}$, with $E$ an open set of real numbers included in $(\alpha,\beta)$, if and only if
$$\overline{\operatorname{Img}\left(f_1(\sigma+it)\right)}=\overline{\operatorname{Img}\left(f_2(\sigma+it)\right)}\ \mbox{for every }\sigma\in(\alpha,\beta).$$
\end{theorem}
\begin{proof}
Suppose that $\overline{\operatorname{Img}\left(f_1(\sigma+it)\right)}=\overline{\operatorname{Img}\left(f_2(\sigma+it)\right)}$
for every $\sigma$ such that $\alpha<\sigma<\beta$. Take an open set $E\subset(\alpha,\beta)$ and $w_0\in \bigcup_{\sigma\in E}\operatorname{Img}\left(f_1(\sigma+it)\right)$, then $w_0\in \operatorname{Img}\left(f_1(\sigma_0+it)\right)$ for some $\sigma_0\in E$ and hence $w_0=f_1(\sigma_0+it_0)$ for some $t_0\in\mathbb{R}$.
Now, by hypothesis, we have $w_0\in\overline{\operatorname{Img}\left(f_2(\sigma_0+it)\right)}$, which yields by Lemma \ref{DeEgea} that the function $f_2(s)$ takes the value $w_0$ on the region $\{s=\sigma+it:\sigma_0-\varepsilon<\sigma<\sigma_0+\varepsilon\}$ for every $\varepsilon>0$ sufficiently small (recall that $E$ is an open set).
Consequently, $w_0\in \bigcup_{\sigma\in E}\operatorname{Img}\left(f_2(\sigma+it)\right)$.
By symmetry, we analogously prove that $\bigcup_{\sigma\in E}\operatorname{Img}\left(f_2(\sigma+it)\right)\subset \bigcup_{\sigma\in E}\operatorname{Img}\left(f_1(\sigma+it)\right)$.
\smallskip

\noindent Conversely, suppose that the functions $f_1$ and $f_2$ take the same set of values on every region $\{s=\sigma+it\in\mathbb{C}:\sigma\in E\}$, where $E$ is an open set in $(\alpha,\beta)$. By reductio ad absurdum, suppose the existence of $\sigma_0\in(\alpha,\beta)$ such that $\overline{\operatorname{Img}\left(f_1(\sigma_0+it)\right)}\neq\overline{\operatorname{Img}\left(f_2(\sigma_0+it)\right)}$.
Thus, without loss of generality, there exists $w_1\in\mathbb{C}$ such that $w_1\in \overline{\operatorname{Img}\left(f_1(\sigma_0+it)\right)}$ and $w_1\notin \overline{\operatorname{Img}\left(f_2(\sigma_0+it)\right)}$. In view of Lemma \ref{DeEgea}, this yields the existence of $\varepsilon_0>0$ such that
$$w_1\in\bigcup_{\sigma\in E_{\sigma_0,\varepsilon}}\operatorname{Img}\left(f_1(\sigma+it)\right)\ \mbox{for every } 0<\varepsilon<\varepsilon_0,$$
where $E_{\sigma_0,\varepsilon}=(\sigma_0-\varepsilon,\sigma_0+\varepsilon)$. Furthermore, since $w_1\notin \overline{\operatorname{Img}\left(f_2(\sigma_0+it)\right)}$, we deduce from the converse of Lemma \ref{DeEgea} the existence of $\varepsilon_1>0$ such that $w_1\notin \bigcup_{\sigma\in E_{\sigma_0,\varepsilon_1}}\operatorname{Img}\left(f_2(\sigma+it)\right)$.
Consequently, by taking $\varepsilon_2=\min\{\varepsilon_0,\varepsilon_1\}$, we conclude that
$$\bigcup_{\sigma\in E_1}\operatorname{Img}\left(f_1(\sigma+it)\right)\neq \bigcup_{\sigma\in E_1}\operatorname{Img}\left(f_2(\sigma+it)\right),$$
where $E_1=(\sigma_0-\varepsilon,\sigma_0+\varepsilon)$ and $0<\varepsilon<\varepsilon_2$. This represents a contradiction and the result follows.
\end{proof}

\begin{corollary}\label{lemma}
Let $f_1\in AP(U_1,\mathbb{C})$ and $f_2\in AP(U_2,\mathbb{C})$, with $U_1=\{\sigma+it\in\mathbb{C}:\alpha_1<\sigma<\beta_1\}$ and $U_2=\{\sigma+it\in\mathbb{C}:\alpha_2<\sigma<\beta_2\}$ such that $U_1\cap U_2\neq \emptyset$. Consider an interval $(\alpha,\beta)\subset(\alpha_1,\beta_1)\cap (\alpha_2,\beta_2)$, and $w\in \mathbb{C}$ with $\operatorname{Re}w\in (\alpha,\beta)$.
If $\overline{\operatorname{Img}\left(f_1(\sigma+it)\right)}=\overline{\operatorname{Img}\left(f_2(\sigma+it)\right)}\ \mbox{for every }\sigma\in(\alpha,\beta)$, then there exist $\{w_m\}_{m\geq 1}\subset U_1\cap U_2$, with $w_m\rightarrow w$, and $\{t_m\}_{m\geq 1}\subset\mathbb{R}$ such that $f_1(w_m+it_m)=f_2(w)$ for each $m\geq 1$.
\end{corollary}
\begin{proof}
Given $w\in \mathbb{C}$ with $\operatorname{Re}w\in (\alpha,\beta)$ and $m\in\mathbb{N}$, consider the set
$$E_m=(\operatorname{Re}w-1/m,\operatorname{Re}w+1/m)\cap (\alpha,\beta).$$
Since $\overline{\operatorname{Img}\left(f_1(\sigma+it)\right)}=\overline{\operatorname{Img}\left(f_2(\sigma+it)\right)}\ \mbox{for every }\sigma\in(\alpha,\beta)$, Theorem \ref{messi} assures that
$$\bigcup_{\sigma\in E}\operatorname{Img}\left(f_1(\sigma+it)\right)=\bigcup_{\sigma\in E}\operatorname{Img}\left(f_2(\sigma+it)\right)$$
for every open subset $E$ in $(\alpha,\beta)$, and in particular for $E_m$ with $m\in\mathbb{N}$.
This yields the existence of at least one point $z_m\in\{s\in\mathbb{C}:\operatorname{Re}s\in E_m\}$ such that $f_2(w)=f_1(z_m)$. Now, if we take $w_m:=\operatorname{Re}z_m+i\operatorname{Im}w$ and $t_m:=\operatorname{Im}z_m-\operatorname{Im}w$, then we have
$$|w-w_m|=|\operatorname{Re}w-\operatorname{Re}z_m|<1/m,$$ so $w_m\rightarrow w$, and $f_2(w)=f_1(w_m+it_m)$.
\end{proof}

\section{On the converse of Bohr's equivalence theorem}

In this section, we will prove a converse of Bohr's equivalence theorem for the case that the Fourier exponents are $\mathbb{Q}$-linearly independent (subsection \ref{sub1}) and for the case that $0$ is a Fourier exponent and the remaining exponents are $\mathbb{Q}$-linearly independent (subsection \ref{sub2}).

\smallskip

Recall that $\sharp \Lambda$ denotes the cardinal of the numerable set $\Lambda$.

\subsection{Sets of exponents linearly independent over the rational numbers}\label{sub1}
\medskip

Given $\Lambda=\{\lambda_1,\lambda_2,\ldots,\lambda_j,\ldots\}$ a set of real numbers which are linearly independent over the rational numbers, consider an open vertical strip of the type $U=\{s\in\mathbb{C}: \alpha<\operatorname{Re}s<\beta\}$, with $-\infty\leq\alpha<\beta\leq\infty$, and $f(s)$ an almost periodic function in $AP(U,\mathbb{C})$ whose Dirichlet series is
of the form
 \begin{equation}\label{eqqo}
\sum_{j\geq 1}a_je^{\lambda_js},\ a_j\in\mathbb{C},\ \lambda_j\in\Lambda.
\end{equation}
Then $f$ can be associated with an auxiliary function $F_f$ of countably many real variables as follows (see \cite[Definition 5]{Sepulcre-Vidal1} and \cite[Definition 6]{Sepulcre-Vidal2} for a more general definition, without $\mathbb{Q}$-linear independence).

\begin{definition}\label{auxuliaryfunc}
Given $\Lambda=\{\lambda_1,\lambda_2,\ldots,\lambda_j,\ldots\}$ a set of exponents which are $\mathbb{Q}$-linearly independent, let $f(s)$ be an almost periodic function in $AP(U,\mathbb{C})$, with $U=\{s\in\mathbb{C}:\alpha<\operatorname{Re}s<\beta\}$, whose Dirichlet series is of the form \eqref{eqqo}.
We define the auxiliary function  $F_f: (\alpha,\beta)
\times
[0,2\pi)
^{\sharp \Lambda}\rightarrow
\mathbb{C}$
 associated with $f$ as
\begin{equation}\label{2.4}
F_{f}(\sigma,\mathbf{x}):=\sum_{j\geq1}a_j e^{\lambda_j\sigma
}e^{x_ji}\text{, }
\end{equation}%
where $\sigma \in
(\alpha,\beta)
\text{, }\mathbf{x}=(x_1,x_2,\ldots)\in[0,2\pi)^{\sharp \Lambda}$ and the series in (\ref{2.4}) is summed by Bochner-Fej\'{e}r procedure, applied at $t=0$ to the sum
$\sum_{j\geq1}a_j e^{x_ji}e^{\lambda_js}$.
\end{definition}

Note that the Dirichlet series $\sum_{j\geq 1}a_je^{\lambda_j(\sigma+it)}$ associated with $f(s)$  arises from its auxiliary function $F_f(\sigma,\mathbf{x}_t)$ by the special choice of $\mathbf{x}_t=t(\lambda_1,\lambda_2,\ldots,\lambda_j,\ldots)+2\pi (p_{t,1},p_{t,2},\ldots,p_{t,j},\ldots)\in[0,2\pi)^{\sharp \Lambda}$ for $p_{t,j}\in\mathbb{Z}$ such that $t\lambda_j+2\pi p_{t,j}\in[0,2\pi)$ for each $j$.
In fact, every arbitrary choice of $\mathbf{x}\in [0,2\pi)^{\sharp \Lambda}$
leads to a Dirichlet series which is {\newequiv}t to that of $f$. 
Moreover, it is worth noting that the condition of $\mathbb{Q}$-linear independence of the Fourier exponents yields by \cite[p. 154]{Besi} that its Dirichlet series is absolutely convergent.
\smallskip

In connection with the auxiliary function $F_f$, we next establish the following notation.

\begin{definition}\label{image}
Given $\Lambda=\{\lambda_1,\lambda_2,\ldots,\lambda_j,\ldots\}$ a set of exponents which are $\mathbb{Q}$-linearly independent, let $f(s)$ be an almost periodic function in $AP(U,\mathbb{C})$ whose Dirichlet series is of the form \eqref{eqqo}, and $\sigma_0=\operatorname{Re}s_0$ with $s_0\in U$. We define $\operatorname{Img}\left(F_f(\sigma_0,\mathbf{x})\right)$ to be the set of values in the complex plane taken on by the auxiliary function $F_f(\sigma,\mathbf{x})$ when $\sigma=\sigma_0$; that is
$$\operatorname{Img}\left(F_f(\sigma_0,\mathbf{x})\right)=\{s\in\mathbb{C}:\exists \mathbf{x}\in[0,2\pi)^{\sharp \Lambda}\mbox{ such that }s=F_f(\sigma_0,\mathbf{x})\}.$$
\end{definition}

Take $f\in AP(U,\mathbb{C})$, $f_1\ \shortstack{$_{{\fontsize{6}{7}\selectfont *}}$\\$\sim$}\ f$ and $\sigma_0=\operatorname{Re}s_0$ with $s_0\in U$. With the notation above,  it was proved in \cite[Lemma 9 and Propositions 12-13]{Sepulcre-Vidal1} (or, more generally, in \cite[Proposition 4]{Sepulcre-Vidal2}) that $\operatorname{Img}\left(F_f(\sigma_0,\mathbf{x})\right)$ is a closed set and
\begin{equation}\label{eqquati}
\operatorname{Img}\left(F_f(\sigma_0,\mathbf{x})\right)=\bigcup_{f_k \shortstack{$_{{\fontsize{6}{7}\selectfont *}}$\\$\sim$} f}\operatorname{Img}\left(f_k(\sigma_0+it)\right)=\overline{\operatorname{Img}\left(f_1(\sigma_0+it)\right)}. 
\end{equation}
In fact, $\operatorname{Img}\left(F_f(\sigma_0,\mathbf{x})\right)$ is a compact set and, if the Dirichlet series of $f$ is of the form \eqref{eqqo}, we have
\begin{equation}\label{kkkiu}
\left|F_f(\sigma_0,\mathbf{x})\right|\leq \sum_{j\geq1}|a_j| e^{\lambda_j\sigma_0
}\ \mbox{for every }\mathbf{x}\in[0,2\pi)^{\sharp \Lambda}.
\end{equation}
It is clear that this maximum value for the modulus of the points in the set $\operatorname{Img}\left(F_f(\sigma_0,\mathbf{x})\right)$ is attained when all the summands of \eqref{2.4} are aligned. In fact, we can prove the following result.

\begin{lemma}\label{lnew}
Given $\Lambda=\{\lambda_1,\lambda_2,\ldots,\lambda_j,\ldots\}$ a set of exponents which are $\mathbb{Q}$-linearly independent, let $f(s)$ be an almost periodic function in $AP(U,\mathbb{C})$ whose Dirichlet series is of the form \eqref{eqqo}, and $\sigma_0=\operatorname{Re}s_0$ with $s_0\in U$. Then the set $$\{s\in\operatorname{Img}\left(F_f(\sigma_0,\mathbf{x})\right): |s|=\max\{\left|F_{f}(\sigma,\mathbf{x})\right|:\mathbf{x}\in[0,2\pi)^{\sharp \Lambda}\}\}$$ coincides with the circumference of centre the origin and radius $\sum_{j\geq1}|a_j| e^{\lambda_j\sigma_0}$.
\end{lemma}
\begin{proof}
Fixed $\sigma_0=\operatorname{Re}s_0$ with $s_0\in U$, we first note that the choice $\mathbf{x}_0=(x_1,x_2,\ldots,x_j,\ldots)$ with $x_j=2\pi-\operatorname{arg}a_j$, $j=1,2,\ldots$ (where $\operatorname{arg}a_j\in (0,2\pi]$) leads to
$$F_f(\sigma_0,\mathbf{x}_0)=\sum_{j\geq1}a_j e^{\lambda_j\sigma
}e^{x_ji}=\sum_{j\geq1}|a_j|e^{i\operatorname{arg}a_j} e^{\lambda_j\sigma
}e^{-i\operatorname{arg}a_j}=\sum_{j\geq1}|a_j| e^{\lambda_j\sigma_0}.$$
In fact, given $\theta\in[0,2\pi)$, by taking the vector $\mathbf{y}_{\theta}=\mathbf{x}_0+\theta\mathbf{1}-2\pi \mathbf{p}_{\theta}$, with the components of $\mathbf{p}_{\theta}$ in $\{0,1\}$ such that $\mathbf{y}_{\theta}\in [0,2\pi)^{\sharp\Lambda}$, we have
$$F_f(\sigma_0,\mathbf{y}_{\theta})=\sum_{j\geq1}a_j e^{\lambda_j\sigma_0
}e^{-\operatorname{arg}a_ji}e^{\theta i}=e^{i\theta}\sum_{j\geq1}|a_j| e^{\lambda_j\sigma_0},$$
which yields that $e^{i\theta}\sum_{j\geq1}|a_j| e^{\lambda_j\sigma_0}\in \operatorname{Img}\left(F_{f}(\sigma_0,\mathbf{x})\right)$ for every $\theta\in[0,2\pi)$. This shows, jointly with \eqref{kkkiu}, that
\begin{equation*}
\begin{split}
&\{s\in\operatorname{Img}\left(F_f(\sigma_0,\mathbf{x})\right): |s|=\max\{\left|F_{f}(\sigma,\mathbf{x})\right|:\mathbf{x}\in[0,2\pi)^{\sharp \Lambda}\}\}\\
&
=\{s\in\operatorname{Img}\left(F_f(\sigma_0,\mathbf{x})\right): |s|=\sum_{j\geq1}|a_j| e^{\lambda_j\sigma_0}\}\\
&=\{s\in\mathbb{C}:|s|=\sum_{j\geq1}|a_j| e^{\lambda_j\sigma_0}\}.
\end{split}
\end{equation*}
\end{proof}

Under the conditions above, we next prove that two almost periodic functions are {\newequiv}t if and only if they
assume the same set of values on every region $\{s=\sigma+it\in\mathbb{C}:\sigma\in E\}$, where $E$ is an arbitrary open subset of the real projections of an open vertical strip included in their common region of almost periodicity. This also shows a converse of \cite[Theorem 1]{Sepulcre-Vidal2} for exponents $\mathbb{Q}$-linearly independent (non-necessarily equal to each other).

\begin{theorem}[Main result]\label{mth}
Consider $U_1=\{\sigma+it\in\mathbb{C}:\alpha_1<\sigma<\beta_1\}$ and $U_2=\{\sigma+it\in\mathbb{C}:\alpha_2<\sigma<\beta_2\}$ with $U_1\cap U_2\neq \emptyset$.
Given $\Lambda_1=\{\lambda_1,\lambda_2,\ldots,\lambda_j,\ldots\}$ and $\Lambda_2=\{\mu_1,\mu_2,\ldots,\mu_j,\ldots\}$ two sets of exponents which are $\mathbb{Q}$-linearly independent, let $f_1\in AP(U_1,\mathbb{C})$ and $f_2\in AP(U_2,\mathbb{C})$ be two almost periodic functions whose Fourier exponents are $\Lambda_1$ and $\Lambda_2$, respectively.
Then $f_1\ \shortstack{$_{{\fontsize{6}{7}\selectfont *}}$\\$\sim$}\ f_2$ if and only if $$\bigcup_{\sigma\in E}\operatorname{Img}\left(f_1(\sigma+it)\right)=\bigcup_{\sigma\in E}\operatorname{Img}\left(f_2(\sigma+it)\right)$$
for every open set $E$ of real numbers included in a certain interval $(\alpha,\beta)\subset(\alpha_1,\beta_1)\cap (\alpha_2,\beta_2)$.
\end{theorem}
\begin{proof}
Suppose that $f_1\in AP(U_1,\mathbb{C})$ and $f_2\in AP(U_2,\mathbb{C})$ are two almost periodic functions whose Dirichlet series are of the form $\sum_{j\geq 1}a_je^{\lambda_js}$ and $\sum_{j\geq 1}b_je^{\mu_js}$, respectively.
We first note that if $f_1$ and $f_2$ are {\newequiv}t, then $U_1=U_2$ and their sets of Fourier exponents are the same. Hence, by \cite[Theorem 1]{Sepulcre-Vidal2} we get $$\bigcup_{\sigma\in E}\operatorname{Img}\left(f_1(\sigma+it)\right)=\bigcup_{\sigma\in E}\operatorname{Img}\left(f_2(\sigma+it)\right)$$
for every open set $E$ of real numbers included in $(\alpha_1,\beta_1)=(\alpha_2,\beta_2)$.
\smallskip

\noindent Conversely, suppose that $f_1$ and $f_2$ take the same set of values on every region $\{s=\sigma+it\in\mathbb{C}:\sigma\in E\}$ for every open set $E$ of real numbers included in a certain interval $(\alpha,\beta)\subset(\alpha_1,\beta_1)\cap (\alpha_2,\beta_2)$. By Theorem \ref{messi}, we get $$\overline{\operatorname{Img}\left(f_1(\sigma+it)\right)}=\overline{\operatorname{Img}\left(f_2(\sigma+it)\right)}\ \mbox{for every }\sigma\in(\alpha,\beta).$$
By \eqref{eqquati}, this means that
$$\operatorname{Img}\left(F_{f_1}(\sigma,\mathbf{x})\right)=\operatorname{Img}\left(F_{f_2}(\sigma,\mathbf{x})\right)\ \mbox{for every }\sigma\in(\alpha,\beta).$$
In particular, $\mbox{for every }\sigma\in(\alpha,\beta)$, we have
$$\max\{\left|F_{f_1}(\sigma,\mathbf{x})\right|:\mathbf{x}\in[0,2\pi)^{\sharp \Lambda}\}=\max\{\left|F_{f_2}(\sigma,\mathbf{x})\right|:\mathbf{x}\in[0,2\pi)^{\sharp \Lambda}\},$$
or, equivalently (see also Lemma \ref{lnew}),
\begin{equation}\label{unodostres}
\sum_{j\geq1}|a_j| e^{\lambda_j\sigma
}=\sum_{j\geq1}|b_j| e^{\mu_j\sigma
}\ \mbox{for every }\sigma\in(\alpha,\beta).
\end{equation}
Now, by \cite[Lemma 3]{Sepulcre-Vidal0}, let $\hat{f}_1(s)$ and $\hat{f}_2(s)$ be the respective almost periodic functions in $AP(U_1,C)$ and $AP(U_2,C)$ whose Dirichlet series are $\sum_{j\geq1}|a_j| e^{\lambda_js}$ and $\sum_{j\geq1}|b_j| e^{\mu_js}$, and $\hat{f}_1\ \shortstack{$_{{\fontsize{6}{7}\selectfont *}}$\\$\sim$}\ f_1$ and $\hat{f}_2\ \shortstack{$_{{\fontsize{6}{7}\selectfont *}}$\\$\sim$}\ f_2$ (as in the proof of Lemma \ref{lnew}, recall that every arbitrary choice of $\mathbf{x}\in [0,2\pi)^{\sharp \Lambda}$ in $F_{f_j}(\sigma,\mathbf{x})$
leads to a Dirichlet series which is {\newequiv}t to that of $f_j$, for $j=1,2$).
Since these Dirichlet series are absolutely convergent (see \cite[pp. 51-52]{Besi} or \cite[p. 154]{Besi}), they are also uniformly convergent and the functions $\hat{f}_1(s)$ and $\hat{f}_2(s)$ coincide with their respective Dirichlet series \cite[p. 144]{Besi}. Consequently, since they are holomorphic in their respective domains, the equality \eqref{unodostres} and the identity principle yield
$$\sum_{j\geq1}|a_j| e^{\lambda_js}=\sum_{j\geq1}|b_j| e^{\mu_js}\ \mbox{for every }s\in \mathbb{C}\ \mbox{with }\operatorname{Re}s\in (\alpha,\beta).$$ In fact, by the uniqueness theorem \cite[p. 148]{Besi}, the functions $\hat{f}_1$ and $\hat{f}_2$ are identical, and the sets $\Lambda_1$ and $\Lambda_2$ of Fourier exponents are equal (and $U_1=U_2$). Consequently,  $f_1$ and $f_2$ are {\newequiv}t.
\end{proof}

\medskip

Now, we can immediately deduce from our main theorem the following particular result for general Dirichlet series (compare with \cite[Theorem C']{Righetti}).

\begin{corollary}\label{gendirser}
Given $\Lambda$ a set of exponents which is $\mathbb{Q}$-linearly independent,
let $f_1(s)$ and $f_2(s)$ be two general Dirichlet series with the same set of Fourier exponents $\Lambda$ and uniformly convergent on the half-plane $\{s=\sigma+it\in\mathbb{C}:\sigma>\alpha\}$ for some real number $\alpha$. Suppose that $f_1(s)$ and $f_2(s)$ take the same set of values on every vertical strip $\{s=\sigma+it\in\mathbb{C}:\alpha<\sigma_0<\sigma<\sigma_1\}$, with $\sigma_0<\sigma_1\leq\infty$.
Then $f_1(s)$ is {\newequiv}t to $f_2(s)$.
\end{corollary}

\medskip

If the Fourier exponents are $\mathbb{Q}$-linearly independent, it is clear that they form an integral basis (see the Preliminaries section). In this case, Bohr-equivalence and \newequiv ce coincide and our main result (Theorem \ref{mth}) can be also formulated in terms of Bohr-equivalent almost periodic functions.

\begin{theorem}[\small{Bohr equivalence theorem and its converse for $\mathbb{Q}$-linearly independent exponents}]\label{th12}
Consider $U_1=\{\sigma+it\in\mathbb{C}:\alpha_1<\sigma<\beta_1\}$ and $U_2=\{\sigma+it\in\mathbb{C}:\alpha_2<\sigma<\beta_2\}$ with $U_1\cap U_2\neq \emptyset$.
Given $\Lambda_1=\{\lambda_1,\lambda_2,\ldots,\lambda_j,\ldots\}$ and $\Lambda_2=\{\mu_1,\mu_2,\ldots,\mu_j,\ldots\}$ two sets of exponents which are $\mathbb{Q}$-linearly independent, let $f_1\in AP(U_1,\mathbb{C})$ and $f_2\in AP(U_2,\mathbb{C})$ be two almost periodic functions whose Fourier exponents are $\Lambda_1$ and $\Lambda_2$, respectively.
Then $f_1$ and $f_2$ are Bohr-equivalent if and only if $$\bigcup_{\sigma\in E}\operatorname{Img}\left(f_1(\sigma+it)\right)=\bigcup_{\sigma\in E}\operatorname{Img}\left(f_2(\sigma+it)\right)$$
for every open set $E$ of real numbers included in a certain interval $(\alpha,\beta)\subset(\alpha_1,\beta_1)\cap (\alpha_2,\beta_2)$.
\end{theorem}

\subsection{Set of exponents of the form $\{0\}\cup \Lambda$, with $\Lambda$ linearly independent over the rational numbers}\label{sub2}
\medskip

We next consider the case that $0$ is a Fourier exponent and the remaining exponents are $\mathbb{Q}$-linearly independent.
In this way, given an open vertical strip of the type $U=\{s\in\mathbb{C}: \alpha<\operatorname{Re}s<\beta\}$, with $-\infty\leq\alpha<\beta\leq\infty$, let $f(s)$ be an almost periodic function in $AP(U,\mathbb{C})$ whose Dirichlet series is
of the form
 \begin{equation}\label{eqqo2}
a_0+\sum_{j\geq 1}a_je^{\lambda_js},\ a_j\in\mathbb{C}\setminus\{0\}\mbox{ for each }j=0,1,2,\ldots,
\end{equation}
where the exponents $\{\lambda_1,\lambda_2,\ldots,\lambda_j,\ldots\}$ are $\mathbb{Q}$-linearly independent.
Then its associated auxiliary function (analogous to that of Definition \ref{auxuliaryfunc}) is defined in the following terms (for a more general case, see \cite[Definition 6]{Sepulcre-Vidal2}).

\begin{definition}\label{auxuliaryfunc2}
Given $\{\lambda_1,\lambda_2,\ldots,\lambda_j,\ldots\}$ a set of exponents which are $\mathbb{Q}$-linearly independent, let $f(s)$ be an almost periodic function in $AP(U,\mathbb{C})$, with $U=\{s\in\mathbb{C}:\alpha<\operatorname{Re}s<\beta\}$, whose Dirichlet series is of the form \eqref{eqqo2}.
We define the auxiliary function  $F_f: (\alpha,\beta)
\times
[0,2\pi)
^{\sharp \Lambda}\rightarrow
\mathbb{C}$
 associated with $f$ as
\begin{equation}\label{2.42}
F_{f}(\sigma,\mathbf{x}):=a_0+\sum_{j\geq1}a_j e^{\lambda_j\sigma
}e^{x_ji}\text{, }
\end{equation}%
where $\sigma \in(\alpha,\beta)
\text{, }\mathbf{x}=(x_1,x_2,\ldots)\in[0,2\pi)^{\sharp \Lambda}$ and the series in (\ref{2.4}) is summed by Bochner-Fej\'{e}r procedure, applied at $t=0$ to the sum
$\sum_{j\geq1}a_j e^{x_ji}e^{\lambda_js}$.
\end{definition}

As in the previous case, every arbitrary choice of the vector $\mathbf{x}\in [0,2\pi)^{\sharp \Lambda}$
leads to a Dirichlet series which is {\newequiv}t to that of $f$. Moreover, the set of values in the complex plane taken on by the auxiliary function $F_f(\sigma,\mathbf{x})$ when $\sigma=\sigma_0\in (\alpha,\beta)$ is defined in the same manner as
$$\operatorname{Img}\left(F_f(\sigma_0,\mathbf{x})\right)=\{s\in\mathbb{C}:\exists \mathbf{x}\in[0,2\pi)^{\sharp \Lambda}\mbox{ such that }s=F_f(\sigma_0,\mathbf{x})\}.$$
If we take $G_{f}(\sigma,\mathbf{x}):=\sum_{j\geq1}a_j e^{\lambda_j\sigma}e^{x_ji}$, where $\sigma \in(\alpha,\beta)$ and $\mathbf{x}=(x_1,x_2,\ldots)\in[0,2\pi)^{\sharp \Lambda}$, then
\begin{equation}\label{auxconan}
\operatorname{Img}\left(F_f(\sigma_0,\mathbf{x})\right)=\{a_0\}+\operatorname{Img}\left(G_f(\sigma_0,\mathbf{x})\right),
\end{equation}
where $\operatorname{Img}\left(G_f(\sigma_0,\mathbf{x})\right)=\{s\in\mathbb{C}:\exists \mathbf{x}\in[0,2\pi)^{\sharp \Lambda}\mbox{ such that }s=G_f(\sigma_0,\mathbf{x})\}.$
That is, the geometric object $\operatorname{Img}\left(F_f(\sigma_0,\mathbf{x})\right)$ is a translation of $\operatorname{Img}\left(G_f(\sigma_0,\mathbf{x})\right)$ with  translation vector given by $a_0$.

Also, if $f\in AP(U,\mathbb{C})$, $f_1\ \shortstack{$_{{\fontsize{6}{7}\selectfont *}}$\\$\sim$}\ f$ and $\sigma_0=\operatorname{Re}s_0$ with $s_0\in U$,  it was proved in \cite[Proposition 4]{Sepulcre-Vidal2} that $\operatorname{Img}\left(F_f(\sigma_0,\mathbf{x})\right)$ is a closed set and
\begin{equation}\label{eqquati2}
\operatorname{Img}\left(F_f(\sigma_0,\mathbf{x})\right)=\bigcup_{f_k \shortstack{$_{{\fontsize{6}{7}\selectfont *}}$\\$\sim$} f}\operatorname{Img}\left(f_k(\sigma_0+it)\right)=\overline{\operatorname{Img}\left(f_1(\sigma_0+it)\right)}. 
\end{equation}
In fact, if the Dirichlet series of $f$ is of the form \eqref{eqqo2}, it is accomplished that
$$\left|F_f(\sigma_0,\mathbf{x})\right|\leq |a_0|+\sum_{j\geq1}|a_j| e^{\lambda_j\sigma_0
}\ \mbox{for every }\mathbf{x}\in[0,2\pi)^{\sharp \Lambda}.$$
This maximum value for the modulus of the points in the set $\operatorname{Img}\left(F_f(\sigma_0,\mathbf{x})\right)$ is attained when all the summands of \eqref{2.42} are aligned.

\smallskip

Now, we can prove the following theorem for the case  that $0$ is a Fourier exponent and the remaining exponents are $\mathbb{Q}$-linearly independent.

\begin{theorem}\label{mth2}
Consider $U_1=\{\sigma+it\in\mathbb{C}:\alpha_1<\sigma<\beta_1\}$ and $U_2=\{\sigma+it\in\mathbb{C}:\alpha_2<\sigma<\beta_2\}$ with $U_1\cap U_2\neq \emptyset$.
Given $\Lambda_1=\{\lambda_1,\lambda_2,\ldots,\lambda_j,\ldots\}$ and $\Lambda_2=\{\mu_1,\mu_2,\ldots,\mu_j,\ldots\}$ two sets of exponents which are $\mathbb{Q}$-linearly independent, let $f_1\in AP(U_1,\mathbb{C})$ and $f_2\in AP(U_2,\mathbb{C})$ be two almost periodic functions whose respective Fourier exponents are $\Lambda_1\cup\{0\}$ and $\Lambda_2\cup\{0\}$. 
Then $f_1\ \shortstack{$_{{\fontsize{6}{7}\selectfont *}}$\\$\sim$}\ f_2$ if and only if $$\bigcup_{\sigma\in E}\operatorname{Img}\left(f_1(\sigma+it)\right)=\bigcup_{\sigma\in E}\operatorname{Img}\left(f_2(\sigma+it)\right)$$
for every open set $E$ of real numbers included in a certain interval $(\alpha,\beta)\subset(\alpha_1,\beta_1)\cap (\alpha_2,\beta_2)$.
\end{theorem}
\begin{proof}
Suppose that $f_1\in AP(U_1,\mathbb{C})$ and $f_2\in AP(U_2,\mathbb{C})$ are two almost periodic functions whose respective Dirichlet series are of the form $a_0+\sum_{j\geq 1}a_je^{\lambda_js}$ and $b_0+\sum_{j\geq 1}b_je^{\mu_js}$, where $\{\lambda_1,\lambda_2,\ldots\}$ and $\{\mu_1,\mu_2,\ldots\}$ are both $\mathbb{Q}$-linearly independent and $a_j,b_j\in\mathbb{C}\setminus\{0\}$ for each $j=0,1,2,\ldots$
As in the proof of Theorem \ref{mth}, we first note that if $f_1$ and $f_2$ are {\newequiv}t, then $U_1=U_2$, their sets of Fourier exponents coincide and, by \cite[Theorem 1]{Sepulcre-Vidal2}, we get the equality under consideration.
\smallskip

\noindent Conversely, suppose that the equality $\bigcup_{\sigma\in E}\operatorname{Img}\left(f_1(\sigma+it)\right)=\bigcup_{\sigma\in E}\operatorname{Img}\left(f_2(\sigma+it)\right)$
is satisfied for every open set $E$ of real numbers included in a certain interval $(\alpha,\beta)\subset(\alpha_1,\beta_1)\cap (\alpha_2,\beta_2)$. By Theorem \ref{messi}, we know that $$\overline{\operatorname{Img}\left(f_1(\sigma+it)\right)}=\overline{\operatorname{Img}\left(f_2(\sigma+it)\right)}\ \mbox{for every }\sigma\in(\alpha,\beta).$$
By \eqref{eqquati2}, this means that
$\operatorname{Img}\left(F_{f_1}(\sigma,\mathbf{x})\right)=\operatorname{Img}\left(F_{f_2}(\sigma,\mathbf{x})\right)\ \mbox{for every }\sigma\in(\alpha,\beta)$.
In fact, under the notation $G_{f_1}(\sigma,\mathbf{x})=\sum_{j\geq1}a_j e^{\lambda_j\sigma}e^{x_ji}$ and $G_{f_2}(\sigma,\mathbf{x})=\sum_{j\geq1}b_j e^{\mu_j\sigma}e^{x_ji}$,
 with $\sigma \in(\alpha,\beta)$ and $\mathbf{x}=(x_1,x_2,\ldots)\in[0,2\pi)^{\sharp \Lambda}$, we deduce from \eqref{auxconan} that
$$\{a_0\}+\operatorname{Img}\left(G_{f_1}(\sigma,\mathbf{x})\right)=\{b_0\}+\operatorname{Img}\left(G_{f_2}(\sigma,\mathbf{x})\right)\ \mbox{for every }\sigma\in(\alpha,\beta).$$
By Lemma \ref{lnew}, recall that the circumferences of centre the origin and radii $\sum_{j\geq1}|a_j| e^{\lambda_j\sigma}$ and $\sum_{j\geq1}|b_j| e^{\mu_j\sigma}$ are respectively included in $\operatorname{Img}\left(G_{f_1}(\sigma,\mathbf{x})\right)$ and $\operatorname{Img}\left(G_{f_2}(\sigma,\mathbf{x})\right)$, and these radii represent the respective maximum values of the modulus of the points in the sets $\operatorname{Img}\left(G_{f_1}(\sigma,\mathbf{x})\right)$ and $\operatorname{Img}\left(G_{f_2}(\sigma,\mathbf{x})\right)$. Particularly, this means that the outer boundary of the region $\{a_0\}+\operatorname{Img}\left(G_{f_1}(\sigma,\mathbf{x})\right)$ (which is the circumference of centre $\{a_0\}$ and radius $\sum_{j\geq1}|a_j| e^{\lambda_j\sigma}$) coincides with the outer boundary of the region $\{b_0\}+\operatorname{Img}\left(G_{f_2}(\sigma,\mathbf{x})\right)$ (which is the circumference of centre $\{b_0\}$ and radius $\sum_{j\geq1}|b_j| e^{\mu_j\sigma}$).
Consequently, the two translated sets (and the two translation vectors) must be equal, which means that $a_0=b_0$ and
\begin{equation}\label{knfvqwl}
\operatorname{Img}\left(G_{f_1}(\sigma,\mathbf{x})\right)=\operatorname{Img}\left(G_{f_2}(\sigma,\mathbf{x})\right)\ \mbox{for every }\sigma\in(\alpha,\beta).
\end{equation}
If we take $g_1(s):=f_1(s)-a_0\in AP(U_1,\mathbb{C})$ and $g_2(s):=f_2(s)-b_0=f_2(s)-a_0\in AP(U_2,\mathbb{C})$ (where all their Fourier exponents are $\mathbb{Q}$-linearly independent), equality \eqref{knfvqwl} is equivalent to
$$\operatorname{Img}\left(F_{g_1}(\sigma,\mathbf{x})\right)=\operatorname{Img}\left(F_{g_2}(\sigma,\mathbf{x})\right)\ \mbox{for every }\sigma\in(\alpha,\beta),
$$
where $F_{g_1}(\sigma,\mathbf{x})$ and $F_{g_2}(\sigma,\mathbf{x})$ are the auxiliary functions associated with $g_1$ and $g_2$, respectively.
Now, by \eqref{eqquati}, we have
$$\overline{\operatorname{Img}\left(g_1(\sigma+it)\right)}=\overline{\operatorname{Img}\left(g_2(\sigma+it)\right)}\ \mbox{for every }\sigma\in(\alpha,\beta)$$
and, by Theorem \ref{messi},
$$\bigcup_{\sigma\in E}\operatorname{Img}\left(g_1(\sigma+it)\right)=\bigcup_{\sigma\in E}\operatorname{Img}\left(g_2(\sigma+it)\right)$$
for every open set $E$ of real numbers included in $(\alpha,\beta)$. Therefore, we deduce from Theorem \ref{mth} that $g_1$ and $g_2$ are {\newequiv}t. Hence $f_1$ and $f_2$ are also {\newequiv}t and the result holds.
\end{proof}

\bigskip

As a conjecture, we think that Theorem \ref{mth} is also true without the condition of $\mathbb{Q}$-linear independence of the Fourier exponents.

\begin{conjecture}
Consider $U_1=\{\sigma+it\in\mathbb{C}:\alpha_1<\sigma<\beta_1\}$ and $U_2=\{\sigma+it\in\mathbb{C}:\alpha_2<\sigma<\beta_2\}$ with $U_1\cap U_2\neq \emptyset$.
Given $\Lambda_1=\{\lambda_1,\lambda_2,\ldots,\lambda_j,\ldots\}$ and $\Lambda_2=\{\mu_1,\mu_2,\ldots,\mu_j,\ldots\}$ two sets of exponents, let $f_1\in AP(U_1,\mathbb{C})$ and $f_2\in AP(U_2,\mathbb{C})$ be two almost periodic functions whose Fourier exponents are $\Lambda_1$ and $\Lambda_2$, respectively.
Then $f_1\ \shortstack{$_{{\fontsize{6}{7}\selectfont *}}$\\$\sim$}\ f_2$ if and only if $$\bigcup_{\sigma\in E}\operatorname{Img}\left(f_1(\sigma+it)\right)=\bigcup_{\sigma\in E}\operatorname{Img}\left(f_2(\sigma+it)\right)$$
for every open set $E$ of real numbers included in a certain interval $(\alpha,\beta)\subset(\alpha_1,\beta_1)\cap (\alpha_2,\beta_2)$.
\end{conjecture}

\noindent \textbf{Acknowledgements.} The first author has been partially supported by a CRM-ISM postdoctoral fellowship and by a fellowship ``Ing. Giorgio Schirillo'' from INdAM. The second author's research has been partially supported by MICIU of Spain under project number PGC2018-097960-B-C22.


\end{document}